\documentclass[reqno,twoside,12pt]{amsart}

\hoffset - 1.5cm

\usepackage{amsmath,amssymb}
\usepackage{wasysym}
\usepackage{xcolor}
\usepackage{hyperref}

\newtheorem{definition}{Definition}[section]
\newtheorem{theorem}[definition]{Theorem}
\newtheorem{lemma}[definition]{Lemma}

\newtheorem{corollary}[definition]{Corollary}
\newtheorem{remark}[definition]{Remark}
\newtheorem{example}[definition]{Example}
\title{Wavelet based solutions to the Poisson and the Helmholtz equations on the $n$-dimensional unit sphere}
\author{Ilona Iglewska-Nowak}\address{West Pomeranian University of Technology, Department of Mathematics, al. Piastów 17, PL-70-310 Szczecin, Poland, ORCID 0000--0002--1938--8055 } \email{iiglewskanowak@zut.edu.pl}
\author{Piotr Stefaniak}\address{Nicolaus Copernicus University, Faculty of Mathematics and Computer Science, ul. Chopina  12/18, PL-87-100 Toru\'n, Poland, ORCID 0000--0002--6117--2573 } \email{cstefan@mat.umk.pl}

\begin{document}

\begin{abstract}
We  present a method of solving partial differential equations on the $n$-dimensional unit sphere using methods based on the continuous wavelet transform derived from approximate identities. We give an explicit analytical solution to the Poisson equation and to the Helmholtz equations. For the first one and for some special values of the parameter in the latter one, we derive a closed formula for the generalized Green function.
\end{abstract}


\date{\today}
\maketitle

\section{Introduction}

In the last thirty years, many authors developed wavelet methods for solving differential equations. Already in the 1990s numerical solutions of ODEs \cite{SYYS98,mF99} and PDEs \cite{vP89,sJ92} on the Euclidean space were found. Some further examples of wavelet application to ODEs are presented in \cite{LH14,HMW19}, and to PDEs in \cite{HMV04,MLZJ05,kU09,HK10,cjG14,ZXZ16,BH18}. Unlike in the articles listed so far, the papers \cite{GS04,mF05,MK08} describe methods for numerical solving of PDEs \emph{on the sphere}. Recent years have brought a number of papers in which wavelet methods were involved to solving fractional differential equations \cite{GR14,WMM14,RBAIS20,AAMASD21}. The list is far from being complete, but it is apparent that wavelet based methods for solving PDEs are usually \emph{numerical}. No publication is known to us that presents an analytical solution.

Theories of continuous spherical wavelets have been developed in the last decades, simultaneously to theories of wavelets over the Euclidean space. Iglewska-Nowak has shown in~\cite[Section~5]{IIN15CWT} that there exist only two essentially different continuous wavelet transforms for spherical signals, namely that based on group theory~\cite{AV,AVn} and that derived from approximate identities~\cite{EBCK09,FGS-book,IIN19CWTC,IIN15CWT,IIN18DW}. In the present paper we show that the latter one can be efficiently applied to solving partial differential equations on the sphere. We give an explicit solution to the Poisson and the Helmholtz equations on the unit sphere of arbitrary dimension in form of a convolution with the Green function. The Green function is given as a series but in the case of the Poisson equation and for some values of the parameter~$a$ in the Helmholtz equation its closed form is derived.

The Helmholtz equation $\Delta^{\!\ast} u+au=f$, with $\Delta^{\!\ast}$ being the Laplacian, naturally appears from general conservation laws of physics and can be interpreted as a wave
equation for monochromatic waves (wave equation in the frequency domain). Applying standard separation of variables, one can reduce some time-dependent equations, including heat conduction equation, Schr\"{o}dinger equation, telegraph and other wavetype, or evolutionary, equations, to a Helmholtz equation with $a\in\mathbb{C}$. Here we study the case $a=k^2$ for some $k^2\in\mathbb{R}$, obtaining a closed form of the Green function for the Helmholtz equation, including some negative values of the parameter~$a$, i.e., a purely imaginary $k$. As far as we know, the closed form for the Green function for negative $a$ has not been given in literature before. This specific case is important for its applications, for instance, it appears in the telegraph equation and corresponds to two types of waves: exponentially growing and exponentially decaying. The relativistic counterpart of the Schr\"{o}dinger equation, called the Klein–Gordan equation, describing a free particle with zero spin in the frequency domain, also leads to the Helmholtz equation with purely imaginary $k$.

The present paper seems to be the first attempt to involve wavelet methods to analytical solving of PDEs. The results concerning the Poisson equation are well known for the $2$-dimensional sphere but to our best knowledge a generalization to the case of Helmholtz equation and to higher dimensions has not been published so far.

The paper is organized as follows. Section~\ref{sec:preliminaries} contains basic information about functions and differential operators on the sphere. In Section~\ref{sec:bwt} we introduce the continuous spherical wavelet transform based on approximate identities. Main results yielding a solution to the Helmholtz and the Poisson equations, contained in Theorem~\ref{thm:PDEsolution} and Theorem~\ref{thm:Green_explicit}, are presented in Section~\ref{sec:PDEs} together with a list of closed representations of the Green function for various dimensions and various values of parameter~$a$. The Appendix~\ref{sec:appendix} contains some integral formulae that can be useful for derivation of closed form of the Green function.

\section{Preliminaries}\label{sec:preliminaries}

A square integrable function~$f$ over the $n$-dimensional unit sphere~$\mathcal S^n\subseteq\mathbb R^{n+1}$, $n\geq 2$, with the rotation-invariant measure~$d\sigma$ normalized such that
$$
\Sigma_n=\int_{\mathcal S^n}d\sigma(x)=\frac{2\pi^{(n+1)/2}}{\Gamma\left((n+1)/2\right)},
$$
can be represented as a Fourier series in terms of the hyperspherical harmonics,
\begin{equation}\label{eq:Fs}
f=\sum_{l=0}^\infty\sum_{k\in\mathcal M_{n-1}(l)}a_l^k(f)\,Y_l^k,
\end{equation}
where $\mathcal M_{n-1}(l)$ denotes the set of sequences $k=(k_0,k_1,\ldots,k_{n-1})$ in $\mathbb N_0^{n-1}\times\mathbb Z$ such that $l\geq k_0\geq k_1\geq\ldots\geq\lvert k_{n-1}\rvert$ and $a_l^k(f)$ are the Fourier coefficients of~$f$. The hyperspherical harmonics of degree~$l$ and order~$k$ are given by
\begin{equation}\label{eq:Ylk}
Y_l^k(x)=A_l^k\prod_{\tau=1}^{n-1}C_{k_{\tau-1}-k_\tau}^{\frac{n-\tau}{2}+k_\tau}(\cos\vartheta_\tau)\sin^{k_\tau}\!\vartheta_\tau\cdot e^{\pm ik_{n-1}\varphi}
\end{equation}
for some constants $A_l^k$. Here, $(\vartheta_1,\dots,\vartheta_{n-1},\varphi)$ are the hyperspherical coordinates of~$x\in\mathcal S^n$,
\begin{align*}
x_1&=\cos\vartheta_1,\\
x_2&=\sin\vartheta_1\cos\vartheta_2,\\
\dots\\
x_n&=\sin\vartheta_1\dots\sin\vartheta_{n-2}\sin\vartheta_{n-1}\cos\varphi,\\
x_{n+1}&=\sin\vartheta_1\dots\sin\vartheta_{n-2}\sin\vartheta_{n-1}\sin\varphi,
\end{align*}
and $\mathcal C_\kappa^K$ are the Gegenbauer polynomials of degree~$\kappa$ and order~$K$.

Zonal (rotation-invariant) functions are those depending only on the first hyperspherical coordinate $\vartheta=\vartheta_1$. Unless it leads to misunderstandings, we identify them with functions of~$\vartheta$ or $t=\cos\vartheta$.
A zonal $\mathcal L^1$-function~$f$ has the following Gegenbauer expansion
\begin{equation}\label{eq:Gegenbauer_expansion}
f(t)=\sum_{l=0}^\infty\widehat f(l)\,C_l^\lambda(t),\qquad t=\cos\vartheta,
\end{equation}
where $\widehat f(l)$ are the Gegenbauer coefficients of~$f$ and $\lambda$ is related to the space dimension by
$$
\lambda=\frac{n-1}{2}.
$$
Consequently, for a zonal $\mathcal L^2$-function $f$ one has
$$
\widehat f(l)=A_l^0\cdot a_l^0(f),
$$
compare~\eqref{eq:Fs}, \eqref{eq:Ylk}, and~\eqref{eq:Gegenbauer_expansion}.

For $f,g\in\mathcal L^1(\mathcal S^n)$, $g$ zonal, their convolution $f\ast g$ is defined by
\begin{equation}\label{eq:convolution}
(f\ast g)(x)=\frac{1}{\Sigma_n}\int_{\mathcal S^n}f(y)\,\tau_xg(y)\,d\sigma(y),\qquad\tau_xg(y)=g(x\cdot y),
\end{equation}
and for $f\in\mathcal L^2(\mathcal S^n)$ it is equal to
$$
f\ast g=\sum_{l=0}^\infty\sum_{k\in\mathcal M_{n-1}(l)} \frac{\lambda}{\lambda+l}\,a_l^k(f)\,\widehat g(l)\,Y_l^k.
$$
If $f$ is a zonal function, then
\begin{equation}\label{eq:convolution_series}
\widehat{f\ast g}(l)=\frac{\lambda}{\lambda+l}\,\widehat f(l)\,\widehat g(l).
\end{equation}

The Laplace-Beltrami operator $\Delta^{\!\ast}$ on the sphere is defined by
\begin{align*}
\Delta^{\!\ast}f&(\vartheta_1,\dots,\vartheta_{n-1},\varphi)\\
&=\sum\limits_{k=1}^{n-1}\left(\prod\limits_{j=1}^k\sin\vartheta_j\right)^{-2}(\sin \vartheta_k)^{k+2-n}
   \frac{\partial}{\partial\vartheta_k}
   \left[\sin^{n-k}\vartheta_k\frac{\partial f(\vartheta_1,\dots,\vartheta_{n-1},\varphi)}{\partial\vartheta_k}\right]\\
&+\left(\prod\limits_{j=1}^k\sin\vartheta_j\right)^{-2}
   \frac{\partial^2f(\vartheta_1,\dots,\vartheta_{n-1},\varphi)}{\partial\varphi^2}.
\end{align*}
It is known that the hyperspherical harmonics are the eigenfunctions of $\Delta^{\!\ast}$, i.e.,
\begin{equation}\label{eq:HH_eigenvectors}
\Delta^{\!\ast} Y_l^k=-l(n+l-1)Y_l^k,
\end{equation}
see \cite[Chapter~II, Theorem 4.1]{Shimakura}.
The relation of $\Delta^{\!\ast}$ and the Laplace operator $\Delta$ is given by
\begin{equation}\label{eq:laplacians}
\Delta f = R^{-n}\frac{\partial}{\partial R}\left(R^{n}\frac{\partial f}{\partial R}\right) + \frac{1}{R^2}\Delta^{\!\ast} f,
\end{equation}
where
$R\geq 0$ is the radial distance of $x\in \mathbb{R}^n$ in the hyperspherical coordinates, see \cite[Chapter~II, Proposition 3.3]{Shimakura}.

The Laplace operator is commutative with $SO(n+1)$-rotations $\Upsilon$,
\begin{equation}\label{eq:commutativity}
\Delta\left[f(\Upsilon x)\right]=(\Delta f)(\Upsilon x),
\end{equation}
see \cite[Chapter~IX, Par.~2, Subsec.~4]{Vilenkin}. Consequently, it follows from \eqref{eq:laplacians} that the same holds for the Laplace-Beltrami operator, see also \cite[Chapter~II, formula (3.15)]{Shimakura}.

Since $\mathcal{S}^n$ is a manifold without boundary, the Green second surface identity implies that for $f,g$
of class $\mathcal C^2$ the following holds:
\begin{equation}\label{eq:green}
\int_{\mathcal{S}^n} \Delta^{\!\ast} f(x) \cdot g(x)\, d\sigma(x)=\int_{\mathcal{S}^n}  f(x) \cdot \Delta^{\!\ast}g(x)\, d\sigma(x).
\end{equation}

The scalar product in $\mathcal L^2(\mathcal S^n)$ is antilinear in the first variable,
$$
\left<f,g\right>=\frac{1}{\Sigma_n}\int_{\mathcal{S}^n}\overline{f(x)}\,g(x)\,d\sigma(x).
$$
Since~$\Delta^{\!\ast}$ is a linear operator, one has
$$
\overline{\Delta^{\!\ast} f}=\Delta^{\!\ast}\overline{f}
$$
and~\eqref{eq:green} can be also written as
\begin{equation}\label{eq:Green_sp}
\left<\Delta^{\!\ast}f,g\right>=\left<f,\Delta^{\!\ast}g\right>.
\end{equation}

\section{The wavelet transform}\label{sec:bwt}

The wavelet transform based on approximate identities has been developed starting from the 1990s by Freeden \emph{et al.} \cite{FGS-book,FW-C,FW} and Bernstein \emph{et al.} \cite{sB09,EBCK09}, as well as by Iglewska-Nowak in the recent years \cite{IIN15CWT,IIN19CWTC}. Its modified version can be found in~\cite{IIN18DW,IIN17FDW}, where two distinct wavelet families are used to analysis and synthesis of a signal. In the present paper we shall use the definition from~\cite{IIN18DW} adapted to the case of rotation-invariant functions with $0$-mean and weight function $\alpha(\rho)=\frac1\rho$ and simplified as in~\cite{IIN19CWTC} (the convergence proof of \cite[Theorem~3.2]{IIN19CWTC} can be easily adapted to the considered case).

\begin{definition}\label{def:bilinear_wavelets} Families $\{\Psi_\rho\}_{\rho\in\mathbb R_+}\subseteq\mathcal L^2(\mathbb S^n)$ and $\{\Omega_\rho\}_{\rho\in\mathbb R_+}\subseteq\mathcal L^2(\mathbb S^n)$ of zonal functions are called an admissible wavelet pair if they satisfy the condition
\begin{equation}\label{eq:wp}\begin{split}
\int_0^\infty\overline{\widehat{\Psi_\rho}(0)}\,\widehat{\Omega_\rho}(0)\,\frac{d\rho}{\rho}&=0,\\
\int_0^\infty\overline{\widehat{\Psi_\rho}(l)}\,\widehat{\Omega_\rho}(l)\,\frac{d\rho}{\rho}&=\left(\frac{\lambda+l}{\lambda}\right)^2,\qquad l\in\mathbb N.
\end{split}\end{equation}
$\{\Omega_\rho\}_{\rho\in\mathbb R_+}$ is the reconstruction wavelet to the wavelet $\{\Psi_\rho\}_{\rho\in\mathbb R_+}$.
\end{definition}

\begin{definition}\label{def:wt} Let $\{\Psi_\rho\}_{\rho\in\mathbb R_+}$ and $\{\Omega_\rho\}_{\rho\in\mathbb R_+}$ be an admissible wavelet pair. Then, the spherical wavelet transform
$$
\mathcal W_\Psi\colon\mathcal L^2(\mathcal S^n)\to\mathcal L^2(\mathbb R_+\times\mathcal S^n)
$$
is defined by
$$
\mathcal W_\Psi f(\rho,y):=\left(f\ast\overline{\Psi_\rho}\right)(y).
$$
\end{definition}

The wavelet transform is invertible (in $\mathcal L^2$-sense) by
$$
f(x)=(\mathcal W_\Omega^{-1}\mathcal W_\Psi f)(x)
   :=\lim_{R\to0}\int_R^{1/R}\left(\mathcal W_\Psi f(\rho,\circ)\ast\Omega_\rho\right)(x)\,\frac{d\rho}{\rho}
$$
for functions~$f\in\mathcal L^2$ with $0$-mean, $\int_{\mathcal S^n}f=0$.\\

If $\Psi=\Omega$, this definition reduces to the classical one~\cite{FGS-book,IIN15CWT}.

\begin{remark} Suppose that the family $\{\Psi_\rho\}_{\rho\in\mathbb R_+}$ satisfies condition
$$
\alpha_l(\Psi):=\left(\frac{\lambda}{\lambda+l}\right)^2\int_0^\infty\lvert\widehat{\Psi_\rho}(l)\rvert^2\,\frac{d\rho}{\rho}\ne0\qquad\text{for }l\in\mathbb N.
$$
Then, it is a wavelet family with the reconstruction wavelet $\{\Omega_\rho\}_{\rho\in\mathbb R_+}$ given by
\begin{equation}\label{eq:rec_wv}\begin{split}
\widehat{\Omega_\rho}(0)&=0,\\
\widehat{\Omega_\rho}(l)&=\frac{1}{\alpha_l(\Psi)}\,\widehat{\Psi_\rho}(l),\quad l\in\mathbb N.
\end{split}\end{equation}
\end{remark}

\begin{example}
An example of wavelets being their own reconstruction family are the Poisson wavelets of order $d\in\mathbb N$~\cite{IIN15PW},
$$
\widehat{g_\rho^d}(l)=\frac{2^d}{\sqrt{\Gamma(2d)}}\cdot(\rho l)^d\,e^{-\rho l}\cdot\frac{\lambda+l}{\lambda},\quad l\in\mathbb N_0.
$$
They are given explicitly as derivatives of the Poisson kernel
\begin{equation}\label{eq:Poisson_kernel}\begin{split}
p_r(y)&=\frac{1}{\Sigma_n}\sum_{l=0}^\infty r^l\,\frac{\lambda+l}{\lambda}\,\mathcal C_l^\lambda(\cos\theta)\\
&=\frac{1}{\Sigma_n}\frac{1-r^2}{(1-2r\cos\theta+r^2)^{(n+1)/2}},\quad r=e^{-\rho},\,y_1=\cos\theta,
\end{split}\end{equation}
by
\begin{equation*}
g_\rho^d=\frac{2^d}{\sqrt{\Gamma(2d)}}\cdot(\rho r\partial_r)^d\,(\Sigma_np_r).
\end{equation*}
\end{example}

\section{The Poisson and the Helmholtz equation}\label{sec:PDEs}

Our aim in this section is to study the Helmholtz equation $\Delta^{\!\ast} u+au=f$ for $a\in\mathbb R$ and the Poisson equation $\Delta^{\!\ast} u=f$.
We start with applying wavelet methods and in the 
second part of the section we give explicit solutions for these equations in form of convolutions with the Green functions.

\subsection{Solution with wavelet methods}

Using the wavelet transform and the inverse wavelet transform with suitably chosen wavelets, we may derive a solution to the Helmholtz and the Poisson equation.

\begin{theorem}\label{thm:PDEsolution}
1) Suppose that $f\in\mathcal C(\mathcal S^n)$ and $u\in\mathcal C^2(\mathcal S^n)$ satisfy
\begin{equation}\label{eq:Helmholtz}
\Delta^{\!\ast} u+au=f,
\end{equation}
where $a\in\mathbb R\setminus\{l(n+l-1),\,l\in\mathbb N_0\}$. Then,
\begin{equation}\label{eq:solution_u}
u=f\ast G,
\end{equation}
for
\begin{equation}\label{eq:Green_function}
G=\sum_{l=0}^\infty\frac{1}{a-l(n+l-1)}\,\frac{\lambda+l}{\lambda}\,\mathcal C_l^\lambda.
\end{equation}
2) Suppose that $f\in\mathcal C(\mathcal S^n)$, $u\in\mathcal C^2(\mathcal S^n)$ are such that $\int f=\int u=0$ and satisfy
$$
\Delta^{\!\ast} u=f.
$$
Then,
\begin{equation}\label{eq:solution_u}
u=f\ast G,
\end{equation}
for
\begin{equation}\label{eq:Green_function}
G=\sum_{l=1}^\infty\frac{-1}{l(n+l-1)}\,\frac{\lambda+l}{\lambda}\,\mathcal C_l^\lambda.
\end{equation}
\end{theorem}

\begin{remark}\begin{enumerate}\item The obtained result coincides with the one presented in~\cite[Chap.~IV, Par.~5]{sgM64}, concerning the general case of `suitable smooth' functions in a vector space.

\item If $a=L(n+L-1)$ for some $L\in\mathbb N_0$, then~\eqref{eq:Helmholtz} can be solved only if
$$
f\ast\mathcal C_L=0.
$$
The solution is unique if one additionally assumes
$$
u\ast\mathcal C_L=0
$$
and it is given by~\eqref{eq:solution_u} with
\begin{equation}\label{eq:Green_function_linN0}
G=\sum_{\l\in\mathbb N_0,\,l\ne L}\frac{1}{a-l(n+l-1)}\,\frac{\lambda+l}{\lambda}\,\mathcal C_l^\lambda,
\end{equation}
compare \cite[Section~4.1]{FS09} or \cite[Section~4.6]{FG13} for the Poisson equation on the two-dimensional sphere (note that the factor~$\frac{1}{4\pi}$ is included in the convolution definition~\eqref{eq:convolution}), resp. \cite{rSz06,rSz07} as well as \cite[Chap.~IV, Par.~5]{sgM64} for the other cases.

\item A solution to the Poisson equation on the $2$-sphere can also be found in~\cite[Theorem~4.7.27]{vM22}. That derivation involves spherical Sobolev spaces and includes more precise considerations on the convergence series for~$u$.

\item A non-zero $\mathcal C^2$-solution to the homogeneous Helmholtz (with $a\ne0$) equation exists on~$\mathcal S^n$ if and only if
$$a=L(n+L-1),$$
for some $L\in\mathbb N$, compare~\cite[Theorem~7.1]{dE68}.
\end{enumerate}\end{remark}

{\bf Proof.} Suppose that $$\int_{\mathcal S^n}f=\int_{\mathcal S^n}u=0$$
and let $\{\Psi_\rho\}$ be the Poisson wavelet family of order~$1$,
\begin{equation}\label{eq:Psi}
\Psi_\rho=2\sum_{l=0}^\infty \rho l\,e^{-\rho l}\cdot\frac{\lambda+l}{\lambda}\,\mathcal C_l^\lambda.
\end{equation}
The wavelet transform of~$f$ with respect to~$\Psi_\rho$ is equal to
$$
\mathcal W_{\Psi} f(\rho,y)=\left<\tau_y\Psi_\rho,f\right>=\left<\tau_y\Psi_\rho,\Delta^{\!\ast} u\right>+\left<\tau_y\Psi_\rho,au\right>,
$$
and it can be written as
\begin{equation}\label{eq:WPsif}
\mathcal W_{\Psi} f(\rho,y)=\left<\tau_y\left[(\Delta^{\!\ast}\!+\bar a)\Psi_\rho\right],u\right>
\end{equation}
(by the commutativity of the Laplace-Beltrami operator with rotations~\eqref{eq:commutativity}). By~\eqref{eq:HH_eigenvectors} and~\eqref{eq:Ylk}, the Gegenbauer coefficients of the family $\Theta_\rho:=(\Delta^{\!\ast}\!+\bar a)\Psi_\rho$ are equal to
\begin{align*}
\widehat{\Theta_\rho}(l)&=[\bar a-l(n+l-1)]\,\widehat{\Psi_\rho}(l)\\
&=2\left[\bar a-l(n+l-1)\right]\cdot\rho le^{-\rho l}\cdot\frac{\lambda+l}{\lambda}
\end{align*}
and they satisfy
\begin{align*}
\int_0^\infty\lvert\widehat{\Theta_\rho}(0)\rvert^2\,\frac{d\rho}{\rho}&=0,\\
\int_0^\infty\lvert\widehat{\Theta_\rho}(l)\rvert^2\,\frac{d\rho}{\rho}
   &=\lvert\bar a-l(n+l-1)\rvert^2\cdot\left(\frac{\lambda+l}{\lambda}\right)^2,\qquad l\in\mathbb N.
\end{align*}
Since $a\ne(n-1)l+l^2$ for all $l\in\mathbb N$, the family $\{\Theta_\rho\}$ is a wavelet with the reconstruction wavelet~$\{\Omega_\rho\}$ given by
\begin{align*}
\widehat{\Omega_\rho}(l)&=\frac{2\rho le^{-\rho l}}{[a-l(n+l-1)]}\cdot\frac{\lambda+l}{\lambda},\quad l\in\mathbb N_0,
\end{align*}
compare~\eqref{eq:rec_wv}.

It follows from~\eqref{eq:WPsif} that
$$
\mathcal W_{\Psi}f=\mathcal W_{\Theta}u
$$
and hence, in the case $\int u=0$,
\begin{align*}
u(x)&=\mathcal W_{\Omega}^{-1}\mathcal W_{\Psi}f(x)\\
&=\int_0^\infty\frac{1}{\Sigma_n}\int_{\mathcal S^n}\Omega_\rho(x\cdot y)
   \cdot\frac{1}{\Sigma_n}\int_{\mathcal S^n}\overline{\Psi_\rho(y\cdot z)}\,f(z)\,d\sigma(z)\,d\sigma(y)\,\frac{d\rho}{\rho}
\end{align*}
with $\int_0^\infty$ understood as $\lim_{R\to0}\int_R^{1/R}$. $\mathcal L^2$-convergence of the triple integral is ensured by the properties of the wavelet families, thus, the order of integration can be changed and one obtains
$$
u=f\ast G,
$$
where
$$
G=\int_0^\infty\left(\overline{\Psi_\rho}\ast\Omega_\rho\right)\frac{d\rho}{\rho}.
$$
Now, by~\eqref{eq:convolution_series},
\begin{align*}
\widehat{\overline{\Psi_\rho}\ast\Omega_\rho}(l)&=2\rho le^{-\rho l}\cdot\frac{2\rho le^{-\rho l}}{[a-l(n+l-1)]}\cdot\frac{\lambda+l}{\lambda}
   =\frac{4(\rho l)^2e^{-2\rho l}}{[a-l(n+l-1)]}\cdot\frac{\lambda+l}{\lambda}.
\end{align*}
Since $\overline{\Psi_\rho}\ast\Omega_\rho\in\mathcal L^2(\mathcal S^n)\subset\mathcal L^1(\mathcal S^n)$, it can be integrated term by term in its Gegenbauer representation, i.e.,
\begin{align*}
\widehat{G}(0)&=0\notag\\
\widehat{G}(l)&=\int_0^\infty\frac{4(\rho l)^2\,e^{-2\rho l}}{[a-l(n+l-1)]}\,\frac{d\rho}{\rho}\cdot\frac{\lambda+l}{\lambda}
=\frac{1}{a-l(n+l-1)}\cdot\frac{\lambda+l}{\lambda}.
\end{align*}
If $\int u\ne0$ and $a\ne0$, it follows directly from~\eqref{eq:Helmholtz} that $\widehat u(0)=\frac{1}{a}\cdot\widehat f(0)$, i.e., $u$ given by~\eqref{eq:solution_u} with~$G$ satisfying~\eqref{eq:Green_function} for all~$l\in\mathbb N_0$ solves the Helmholtz equation.
\hfill$\Box$\\

\subsection{Closed forms of the Green functions}

In the papers \cite{rSz06,rSz07} closed forms of the Green functions~\eqref{eq:Green_function_linN0} for the Helmholtz equation with $a=L(n+L-1)$ for some $L\in\mathbb N_0$ are given. In the remaining part of the paper we develop a method for finding closed expressions of these functions. It is based on the series representation of the Poisson kernel and has been used for finding explicit formulae of various wavelets. Contrary to the one from \cite{rSz06,rSz07}, it is applicable also in the case $a\ne L(n+L-1)$. For the set of indices that we have tested, the expressions for~$G$ coincide with those derived by Szmytkowski. For the Poisson equation on the $2$-sphere we obtain the same Green function as the authors of \cite{FS09,FG13,vM22}.

\begin{theorem}\label{thm:Green_explicit} Let $n\in\mathbb N$, $n\geq2$, be fixed, $\lambda=\frac{n-1}{2}$, and suppose that $a=L(n+L-1)$ for some $L\in\mathbb R\setminus\mathbb Z$. Let  $L_0:=\max\{\left[L\right],\left[-n-L+1\right]\}$, where $\left[x\right]$ stays for the biggest integer less than or equal to~$x$. Denote by~$G$ the function
\begin{equation}\label{eq:Green}
G:=\sum_{l=0}^\infty\frac{1}{a-l(n+l-1)}\frac{\lambda+l}{\lambda}\,\mathcal C_l^\lambda.
\end{equation}
Then
\begin{equation}\label{eq:GL}\begin{split}
G(\cos\vartheta)&=-\int_0^1 R^{-(n+2L)}\int_0^Rr^{n+L-2}\\
&\cdot\left(\Sigma_n\cdot p_r(\cos\vartheta)-\sum_{l=0}^{L_0}r^l\cdot\frac{\lambda+l}{\lambda}\,
   \mathcal C_l^\lambda(\cos\vartheta)\right)\,dr\,dR\\
&+\sum_{l=0}^{L_0}\frac{1}{a-l(n+l-1)}\cdot\frac{\lambda+l}{\lambda}\,\mathcal C_l^\lambda(\cos\vartheta).
\end{split}\end{equation}
(If $L_0<0$, set $\sum_0^{L_0}=0$).
\end{theorem}

\begin{bfseries}Proof.\end{bfseries}
According to~\eqref{eq:Poisson_kernel}, the integrand
$$
\Sigma_n\cdot p_r(\cos\vartheta)-\sum_{l=0}^{L_0}r^l\cdot\frac{\lambda+l}{\lambda}\mathcal C_l^\lambda(\cos\vartheta)
$$
equals
\begin{equation}\label{eq:integrand}
\sum_{l=L_0+1}^\infty r^l\cdot\frac{\lambda+l}{\lambda}C_l^\lambda(\cos\vartheta).
\end{equation}
Further, for 
\begin{equation}\label{eq:big_l}
l>\max\{-n-L+1,L\}
\end{equation}
(i.e., $l>L_0$) we have
\begin{equation}\label{eq:int_rR}
\int_0^1R^{-(n+2L)}\int_0^R r^{n+L-2}\cdot r^l dr\,dR=\frac{-1}{(L-l)(n+L+l-1)}.
\end{equation}
Condition~\eqref{eq:big_l} ensures convergence of the integral on the left-hand-side of~\eqref{eq:int_rR} and well definiteness of the quotient 
$$
\frac{-1}{(L-l)(n+L+l-1)}=\frac{-1}{a-l(n+l-1)}.
$$
Since the Gegenbauer polynomials~$\mathcal C_l^\lambda$ over the interval $[-1,1]$ are bounded by
\begin{equation}\label{eq:Cllambda_bound}
\lvert\mathcal C_l^\lambda(\cos\vartheta)\rvert\leq(n+l-2)^{n-2}
\end{equation}
uniformly in $\vartheta\in[0,\pi]$ (compare \cite[Theorem~7.33.1]{Sz75}), the series~\eqref{eq:integrand} is absolutely convergent for $r\in[0,1)$. Consequently, the first summand in~\eqref{eq:GL} equals
\begin{align*}
-\int_0^1&R^{-(n+2L)}\int_0^R r^{n+L-2}\sum_{l=L_0+1}^\infty r^l
   \cdot\frac{\lambda+l}{\lambda}\mathcal C_l^\lambda(\cos\vartheta)\,dr\,dR\\
&=\sum_{l=L_0+1}^\infty\left(-\int_0^1R^{-(n+2L)}\int_0^R r^{n+L-2} r^ldr\,dR\right)
   \cdot\frac{\lambda+l}{\lambda}\mathcal C_l^\lambda(\cos\vartheta)\\
&=\sum_{l=L_0+1}^\infty\frac{1}{a-l(n+l-1)}\cdot\frac{\lambda+l}{\lambda}\,\mathcal C_l^\lambda(\cos\vartheta)
\end{align*}
and~\eqref{eq:GL} coincides with~\eqref{eq:Green}.\hfill$\Box$

\begin{corollary}\label{cor:Green_explicit} Let $n\in\mathbb N$, $n\geq2$, be fixed, $\lambda=\frac{n-1}{2}$, and suppose that $a=L(n+L-1)$ for some $L\in\mathbb N_0$. Further, let~$G$ denote the function
\begin{equation}\label{eq:Green1}
G=\sum_{l=0,l\ne L}^\infty\frac{1}{a-l(n+l-1)}\frac{\lambda+l}{\lambda}\,\mathcal C_l^\lambda
\end{equation}
 Then
\begin{equation*}\begin{split}
G(\cos\vartheta)&=-\int_0^1 R^{-(n+2L)}\int_0^Rr^{n+L-2}\\
&\cdot\left(\Sigma_n\cdot p_r(\cos\vartheta)-\sum_{l=0}^{L}r^l\cdot\frac{\lambda+l}{\lambda}\,
   \mathcal C_l^\lambda(\cos\vartheta)\right)\,dr\,dR\\
&+\sum_{l=0}^{L-1}\frac{1}{a-l(n+l-1)}\cdot\frac{\lambda+l}{\lambda}\,\mathcal C_l^\lambda(\cos\vartheta).
\end{split}\end{equation*}
\end{corollary}

\begin{remark} If $a=L(n+L-1)$, $L\in\mathbb Z$, then also $a=L^\prime(n+L^\prime-1)$ for $L^\prime=-n-L+1$. For nonzero~$a$, the numbers $L$ and $L^\prime$ are of opposite sign and one can choose~$L$ to be positive in order to apply Corollary~\ref{cor:Green_explicit}. If $a=0$, take $L=0$.
\end{remark}

Table~\ref{tab:kernel_K} gives the expressions for~$G$ for $n=2,3,4,5,6,7,8,9,10$.

\begin{table}\caption{The Green function of the Poisson equation~$G$ for different values of~$n$, $t=\cos\vartheta$}\label{tab:kernel_K}
\vspace{0.5em}\centering\begin{tabular}{|l|l|}
\hline
$n=2$&$1+\ln\left(\frac{1-t}{2}\right)$\\[0.5em]
$n=3$&$\frac{-(\pi-\vartheta)\cdot t}{2\sqrt{1-t^2}}+\frac{1}{4}$\\[0.5em]
$n=4$&$\frac{4-7t}{9(1-t)}+\frac{1}{3}\ln\left(\frac{1-t}{2}\right)$\\[0.5em]
$n=5$&$\frac{(\pi-\vartheta)\cdot t(3-2t^2)}{2(1-t^2)^{3/2}}-\frac{3-5t^2}{16(1-t^2)}$\\[0.5em]
$n=6$&$\frac{23-71t+43t^2}{75(1-t)^2}+\frac{1}{5}\ln\left(\frac{1-t}{2}\right)$\\[0.5em]
$n=7$&$\frac{(\pi-\vartheta)\cdot(-15+20t^2-8t^4)}{48(1-t^2)^{5/2}}+\frac{22-71t^2+40t^4}{144(1-t^2)^2}$\\[0.5em]
$n=8$&$\frac{176-759t+906t^2-337t^3}{735(1-t)^3}+\frac{1}{7}\ln\left(\frac{1-t}{2}\right)$\\[0.5em]
$n=9$&$\frac{(\pi-\vartheta)\cdot(-35+70t^2-56t^4+16t^6)}{128(1-t^2)^{7/2}}+\frac{50-237t^2+266t^4-94t^6}{384(1-t^2)^3}$\\[0.5em]
$n=10$&$\frac{563-3089t+5466t^2-4049t^3+1091t^4}{2835(1-t)^4}+\frac{1}{9}\ln\left(\frac{1-t}{2}\right)$\\\hline
\end{tabular}\end{table}

\begin{example} For $n=2$ and $a=0$ (the Poisson equation) we have
$$
\beta_r(t):=\sum_{l=1}^\infty r^{l-1}\cdot(2l+1)\,\mathcal C_l^{1/2}(t)=\frac{1-r^2}{r(1-2tr+r^2)^{3/2}}-\frac{1}{r}.
$$
Further,
$$
\gamma_r(t):=\int\beta_r(t)\,dr=\frac{2}{\sqrt{1-2tr+r^2}}-\ln\left(1-rt+\sqrt{1-2tr+r^2}\right)+C
$$
and
\begin{align*}
\zeta_R(t):=&\int[\gamma_R(t)-\gamma_0(t)]\,dR\\
   =&\ln\left(R-t+\sqrt{1-2tR+R^2}\right)-R\left(1+\ln\frac{1-tR+\sqrt{1-2tR+R^2}}{2}\right).
\end{align*}
Consequently,
$$
G(t)=\zeta_0(t)-\zeta_1(t)=1+\ln\frac{1-t}{2}.
$$
For $t=\cos\vartheta$ it can be expressed as
$$
G(\cos\vartheta)=1+\ln\left(\left(\sin\frac{\vartheta}{2}\right)^2\right).
$$
This coincides (up to the constant $\frac{1}{4\pi}$) with the result from \cite[Lemma~4.3]{FS09} or \cite[Lemma~4.6.2]{FG13}.
\end{example}

Table~\ref{tab:GL_positive_integer} gives the expressions for~$G(\cos\vartheta)$ for several values of parameters $n$ and positive integer $L$. They coincide with those derived in~\cite{rSz06} and~\cite{rSz07}.
\begin{table}\caption{The Green function~$G$ for different values of~$n$ and positive integer $L$; $t=\cos\vartheta$}\label{tab:GL_positive_integer}
\vspace{0.5em}\centering\begin{tabular}{|lll|l|}
\hline
$n=2$, &\hspace{-1em}$L=1$, &\hspace{-1em}$a=2$ &$1+\frac{4}{3}t+t\cdot\ln\left(\frac{1-t}{2}\right)$\\[0.5em]
$n=2$, &\hspace{-1em}$L=2$, &\hspace{-1em}$a=6$ &$\frac{-7+30t+41t^2}{20}-\frac{1-3t^2}{2}\cdot\ln\left(\frac{1-t}{2}\right)$\\[0.5em]
$n=2$, &\hspace{-1em}$L=3$, &\hspace{-1em}$a=12$ &$\frac{-56-123t+210t^2+289t^3}{84}-\frac{t(3-5t^2)}{2}\cdot\ln\left(\frac{1-t}{2}\right)$\\[0.5em]
$n=2$, &\hspace{-1em}$L=4$, &\hspace{-1em}$a=20$ &$\frac{75-660t-1182t^2+1260t^3+1739t^4}{288}+\frac{3-30t^2+35t^4}{8}\cdot\ln\left(\frac{1-t}{2}\right)$\\[0.5em]
$n=3$, &\hspace{-1em}$L=1$, &\hspace{-1em}$a=3$ &$\frac{(\pi-\vartheta)\cdot(1-2t^2)}{2\sqrt{1-t^2}}+\frac{t}{4}$\\[0.5em]
$n=3$, &\hspace{-1em}$L=2$, &\hspace{-1em}$a=8$ &$\frac{(\pi-\vartheta)\cdot t(3-4t^2)}{2\sqrt{1-t^2}}-\frac{1-4t^2}{12}$\\[0.5em]
$n=3$, &\hspace{-1em}$L=3$, &\hspace{-1em}$a=15$ &$\frac{(\pi-\vartheta)\cdot(-1+8t^2-8t^4)}{2\sqrt{1-t^2}}-\frac{t(1-2t^2)}{4}$\\[0.5em]
$n=3$, &\hspace{-1em}$L=4$, &\hspace{-1em}$a=24$ &$\frac{(\pi-\vartheta)\cdot t(-5+20t^2-16t^4)}{2\sqrt{1-t^2}}+\frac{1-12t^2+16t^4}{20}$\\[0.5em]
$n=4$, &\hspace{-1em}$L=1$, &\hspace{-1em}$a=4$ &$\frac{10+13t-28t^2}{15(1-t)}+t\cdot\ln\left(\frac{1-t}{2}\right)$\\[0.5em]
$n=4$, &\hspace{-1em}$L=2$, &\hspace{-1em}$a=10$ &$\frac{-41+223t+149t^2-359t^3}{84(1-t)}-\frac{1-5t^2}{2}\cdot\ln\left(\frac{1-t}{2}\right)$\\[0.5em]
$n=4$, &\hspace{-1em}$L=3$, &\hspace{-1em}$a=18$ &$\frac{-96-213t+903t^2+397t^3-1027t^4}{108(1-t)}-\frac{5t(3-7t^2)}{6}\cdot\ln\left(\frac{1-t}{2}\right)$\\[0.5em]
$n=4$, &\hspace{-1em}$L=4$, &\hspace{-1em}$a=28$ &$\frac{577-5549t-6406t^2+24886t^3+8069t^4-21929t^5}{1056(1-t)}+\frac{5(1-14t^2+21t^4)}{8}\cdot\ln\left(\frac{1-t}{2}\right)$\\[0.5em]
$n=5$, &\hspace{-1em}$L=1$, &\hspace{-1em}$a=5$ &$\frac{(\pi-\vartheta)\cdot(3-12t^2+8t^4)}{8(1-t^2)^{3/2}}+\frac{t(13-16t^2)}{24(1-t^2)}$\\[0.5em]
$n=5$, &\hspace{-1em}$L=2$, &\hspace{-1em}$a=12$ &$\frac{(\pi-\vartheta)\cdot t(15-40t^2+24t^4)}{8(1-t^2)^{3/2}}-\frac{3-23t^2+22t^4}{16(1-t^2)}$\\[0.5em]
$n=5$, &\hspace{-1em}$L=3$, &\hspace{-1em}$a=21$ &$\frac{(\pi-\vartheta)\cdot(-5+60t^2-120t^4+64t^6)}{8(1-t^2)^{3/2}}-\frac{t(37-144t^2+112t^4)}{40(1-t^2)}$\\[0.5em]
$n=5$, &\hspace{-1em}$L=4$, &\hspace{-1em}$a=32$ &$\frac{(\pi-\vartheta)\cdot t(-35+210t^2-336t^4+160t^6)}{8(1-t^2)^{3/2}}+\frac{9-159t^2+416t^4-272t^6}{48(1-t^2)}$\\[0.5em]
$n=6$, &\hspace{-1em}$L=1$, &\hspace{-1em}$a=6$ &$\frac{56+64t-359t^2+232t^3}{105(1-t)^2}+t\cdot\ln\left(\frac{1-t}{2}\right)$\\[0.5em]
$n=6$, &\hspace{-1em}$L=2$, &\hspace{-1em}$a=14$ &$\frac{-103+692t+6t^2-1844t^3+1237t^4}{180(1-t)^2}-2(1-7t^2)\cdot\ln\left(\frac{1-t}{2}\right)$\\[0.5em]
$n=6$, &\hspace{-1em}$L=3$, &\hspace{-1em}$a=24$ &$\frac{-704-1519t+11288t^2-3342t^3-18464t^4+12697t^5}{660(1-t)^2}-\frac{7t(1-3t^2)}{2}\cdot\ln\left(\frac{1-t}{2}\right)$\\[0.5em]
$n=6$, &\hspace{-1em}$L=4$, &\hspace{-1em}$a=36$ &$\frac{5477-58742t-30293t^2+384684t^3-166405t^4-450454t^5+315317t^6}{6240(1-t)^2}$\\[0.5em]
   &&&$+\frac{7(1-18t^2+33t^4)}{8}\cdot\ln\left(\frac{1-t}{2}\right)$\\[0.5em]
$n=7$, &\hspace{-1em}$L=1$, &\hspace{-1em}$a=7$ &$\frac{(\pi-\vartheta)\cdot(5-30t^2+40t^4-16t^6)}{16(1-t^2)^{5/2}}+\frac{t(35-84t^2+46t^4)}{48(1-t^2)^2}$\\[0.5em]
$n=7$, &\hspace{-1em}$L=2$, &\hspace{-1em}$a=16$ &$\frac{(\pi-\vartheta)\cdot(35-140t^2+168t^4-64t^6)}{16(1-t^2)^{5/2}}-\frac{62-695t^2+1304t^4-656t^6}{240(1-t^2)^2}$\\[0.5em]
$n=7$, &\hspace{-1em}$L=3$, &\hspace{-1em}$a=27$ &$\frac{(\pi-\vartheta)\cdot(-35+560t^2-1680t^4+1792t^6-640t^8)}{48(1-t^2)^{5/2}}-\frac{t(255-1462t^2+2240t^4-1024t^6)}{144(1-t^2)^2}$\\[0.5em]
$n=7$, &\hspace{-1em}$L=4$, &\hspace{-1em}$a=40$ &$\frac{(\pi-\vartheta)\cdot t(-105+840t^2-2016t^4+1920t^6-640t^8)}{16(1-t^2)^{5/2}}$\\[0.5em]
   &&&$+\frac{122-2831t^2+10960t^4-14160t^6+5888t^8}{336(1-t^2)^2}$\\[0.5em]
$n=8$, &\hspace{-1em}$L=1$, &\hspace{-1em}$a=8$ &$\frac{144+131t-1518t^2+2013t^3-776t^4}{315(1-t)^3}+t\cdot\ln\left(\frac{1-t}{2}\right)$\\[0.5em]
$n=8$, &\hspace{-1em}$L=2$, &\hspace{-1em}$a=18$ &$\frac{-2927+23367t-14646t^2-75134t^3+114273t^4-45021t^5}{4620(1-t)^3}-\frac{1-9t^2}{2}\cdot\ln\left(\frac{1-t}{2}\right)$\\[0.5em]
$n=8$, &\hspace{-1em}$L=3$, &\hspace{-1em}$a=30$ &$\frac{-6656-13551t+155859t^2-155858t^3-250170t^4+450453t^5-180181t^6}{5460(1-t)^3}$\\[0.5em]
   &&&$-\frac{3t(3-11t^2)}{2}\cdot\ln\left(\frac{1-t}{2}\right)$\\[0.5em]
$n=8$, &\hspace{-1em}$L=4$, &\hspace{-1em}$a=44$ &$\frac{4173-49699t+5793t^2+402945t^3-485545t^4-379929t^5+843387t^6-341189t^7}{3360(1-t)^3}$\\[0.5em]
   &&&$+\frac{3(3-66t^2-143t^4)}{8}\cdot\ln\left(\frac{1-t}{2}\right)$\\[0.5em]\hline
\end{tabular}\end{table}

In Table~\ref{tab:GL_positive_rational} closed forms of the Green function for the Helmholtz equation are given for some positive rational (noninteger) values of the parameter~$a$.
\begin{table}\caption{The Green function~$G$ for different values of~$n$ and positive rational $a$; $t=\cos\vartheta$}\label{tab:GL_positive_rational}
\vspace{0.5em}\centering\begin{tabular}{|lll|l|}
\hline
$n=3$,&\hspace{-1em}$L=\frac{1}{2}$,&\hspace{-1em}$a=\frac{5}{4}$&$\frac{1-2t}{4\sqrt2}\left(\frac{1+t}{\sqrt{1-t}}+\sqrt{1-t}\right)\pi$\\[0.5em]
$n=3$,&\hspace{-1em}$L=\frac{3}{2}$,&\hspace{-1em}$a=\frac{21}{4}$&$\frac{1+2t-4t^2}{4\sqrt2}\left(\frac{1+t}{\sqrt{1-t}}+\sqrt{1-t}\right)\pi$\\[0.5em]
$n=3$,&\hspace{-1em}$L=\frac{5}{2}$,&\hspace{-1em}$a=\frac{45}{4}$&$\frac{-1+4t+4t^2-8t^3}{4\sqrt2}\left(\frac{1+t}{\sqrt{1-t}}+\sqrt{1-t}\right)\pi$\\[0.5em]
$n=3$,&\hspace{-1em}$L=\frac{7}{2}$,&\hspace{-1em}$a=\frac{77}{4}$&$\frac{-1-4t+12t^2+8t^3-16t^4}{4\sqrt2}\left(\frac{1+t}{\sqrt{1-t}}+\sqrt{1-t}\right)\pi$\\[0.5em]
$n=5$,&\hspace{-1em}$L=\frac{1}{2}$,&\hspace{-1em}$a=\frac{9}{4}$&$\frac{-5+18t-12t^2}{32\sqrt2(-1+t)}\left(\frac{1+t}{\sqrt{1-t}}+\sqrt{1-t}\right)\pi$\\[0.5em]
$n=5$,&\hspace{-1em}$L=\frac{3}{2}$,&\hspace{-1em}$a=\frac{33}{4}$&$\frac{-7-12t+60t^2-40t^3}{32\sqrt2(-1+t)}\left(\frac{1+t}{\sqrt{1-t}}+\sqrt{1-t}\right)\pi$\\[0.5em]
$n=5$,&\hspace{-1em}$L=\frac{5}{2}$,&\hspace{-1em}$a=\frac{65}{4}$&$\frac{9-52t-12t^2+168t^3-112t^4}{32\sqrt2(-1+t)}\left(\frac{1+t}{\sqrt{1-t}}+\sqrt{1-t}\right)\pi$\\[0.5em]
$n=5$,&\hspace{-1em}$L=\frac{7}{2}$,&\hspace{-1em}$a=\frac{105}{4}$&$\frac{11+42t-228t^2+32t^3+432t^4-288t^5}{32\sqrt2(-1+t)}\left(\frac{1+t}{\sqrt{1-t}}+\sqrt{1-t}\right)\pi$\\[0.5em]
$n=7$,&\hspace{-1em}$L=\frac{1}{2}$,&\hspace{-1em}$a=\frac{13}{4}$&$\frac{15-76t+100t^2-40t^3}{256\sqrt2(1-t)^2}\left(\frac{1+t}{\sqrt{1-t}}+\sqrt{1-t}\right)\pi$\\[0.5em]
$n=7$,&\hspace{-1em}$L=\frac{3}{2}$,&\hspace{-1em}$a=\frac{45}{4}$&$\frac{77+100t-1020t^2+1400t^3-560t^4}{384\sqrt2(1-t)^2}\left(\frac{1+t}{\sqrt{1-t}}+\sqrt{1-t}\right)\pi$\\[0.5em]
$n=7$,&\hspace{-1em}$L=\frac{5}{2}$,&\hspace{-1em}$a=\frac{85}{4}$&$\frac{-39+290t-140t^2-1120t^3+1680t^4-672t^5}{128\sqrt2(1-t)^2}\left(\frac{1+t}{\sqrt{1-t}}+\sqrt{1-t}\right)\pi$\\[0.5em]
$n=7$,&\hspace{-1em}$L=\frac{7}{2}$,&\hspace{-1em}$a=\frac{133}{4}$&$\frac{-55-194t+1640t^2-1440t^3-3120t^4+5280t^5-2112t^6}{128\sqrt2(1-t)^2}\left(\frac{1+t}{\sqrt{1-t}}+\sqrt{1-t}\right)\pi$\\[0.5em]\hline
\end{tabular}\end{table}

Table~\ref{tab:GL_negative} collects closed forms of the Green function for the Helmholtz equation for some negative values of the parameter~$a$.
\begin{table}\caption{The Green function~$G$ for different values of~$n$ and negative $a$; $t=\cos\vartheta$}\label{tab:GL_negative}
\vspace{0.5em}\centering\begin{tabular}{|lll|l|}
\hline
$n=3$,&\hspace{-1em}$L=-\frac{1}{2}$,&\hspace{-1em}$a=-\frac{3}{4}$&$\frac{-1}{4\sqrt2}\left(\frac{1+t}{\sqrt{1-t}}+\sqrt{1-t}\right)\pi$\\[0.5em]
$n=4$,&\hspace{-1em}$L=-1$,&\hspace{-1em}$a=-2$&$\frac{-1}{3-3t}$\\[0.5em]
$n=5$,&\hspace{-1em}$L=-\frac{1}{2}$,&\hspace{-1em}$a=-\frac{7}{4}$&$\frac{-3+2t}{32\sqrt2(1-t)}\left(\frac{1+t}{\sqrt{1-t}}+\sqrt{1-t}\right)\pi$\\[0.5em]
$n=5$,&\hspace{-1em}$L=-1$,&\hspace{-1em}$a=-3$&$\frac{-t}{8(1-t^2)}-\frac{\vartheta}{8(1-t^2)^{3/2}}$\\[0.5em]
$n=5$,&\hspace{-1em}$L=-\frac{3}{2}$,&\hspace{-1em}$a=-\frac{15}{4}$&$\frac{-1}{32\sqrt2(1-t)}\left(\frac{1+t}{\sqrt{1-t}}+\sqrt{1-t}\right)\pi$\\[0.5em]
$n=5$,&\hspace{-1em}$L=-2$,&\hspace{-1em}$a=-4$&$\frac{-1}{8(1-t^2)}-\frac{t\cdot\vartheta}{8(1-t^2)^{3/2}}$\\[0.5em]
$n=6$,&\hspace{-1em}$L=-1$,&\hspace{-1em}$a=-4$&$\frac{-2+t}{15(1-t)^2}$\\[0.5em]
$n=6$,&\hspace{-1em}$L=-2$,&\hspace{-1em}$a=-6$&$\frac{-1}{15(1-t)^2}$\\[0.5em]
$n=7$,&\hspace{-1em}$L=-\frac{1}{2}$,&\hspace{-1em}$a=-\frac{11}{4}$&$\frac{-7+10t-4t^2}{256\sqrt2(1-t)^2}\left(\frac{1+t}{\sqrt{1-t}}+\sqrt{1-t}\right)\pi$\\[0.5em]
$n=7$,&\hspace{-1em}$L=-1$,&\hspace{-1em}$a=-5$&$\frac{-5t+7t^3-2t^5}{48(1-t^2)^3}-\frac{\vartheta}{16(1-t^2)^{5/2}}$\\[0.5em]
$n=7$,&\hspace{-1em}$L=-\frac{3}{2}$,&\hspace{-1em}$a=-\frac{27}{4}$&$\frac{-5+2t}{768\sqrt2(1-t)^2}\left(\frac{1+t}{\sqrt{1-t}}+\sqrt{1-t}\right)\pi$\\[0.5em]
$n=7$,&\hspace{-1em}$L=-2$,&\hspace{-1em}$a=-8$&$\frac{-2+t^2+t^4}{48(1-t^2)^3}-\frac{t\cdot\vartheta}{16(1-t^2)^{5/2}}$\\[0.5em]
$n=7$,&\hspace{-1em}$L=-\frac{5}{2}$,&\hspace{-1em}$a=-\frac{35}{4}$&$\frac{-1}{128\sqrt2(1-t)^2}\left(\frac{1+t}{\sqrt{1-t}}+\sqrt{1-t}\right)\pi$\\[0.5em]
$n=7$,&\hspace{-1em}$L=-3$,&\hspace{-1em}$a=-9$&$\frac{-t}{16(1-t^2)^2}-\frac{(1+2t^2)\vartheta}{48(1-t^2)^{5/2}}$\\[0.5em]
$n=8$,&\hspace{-1em}$L=-1$,&\hspace{-1em}$a=-6$&$\frac{-8+9t-3t^2}{105(1-t)^3}$\\[0.5em]
$n=8$,&\hspace{-1em}$L=-2$,&\hspace{-1em}$a=-10$&$\frac{-3+t}{105(1-t)^3}$\\[0.5em]
$n=8$,&\hspace{-1em}$L=-3$,&\hspace{-1em}$a=-12$&$\frac{-2}{105(1-t)^3}$\\[0.5em]\hline
\end{tabular}\end{table}\\

Theorem~\ref{thm:Green_explicit} and Corollary~\ref{cor:Green_explicit} hold for all dimensions and all $a$'s, but a closed expression on the right-hand-side of~\eqref{eq:GL} can be rarely obtained for noninteger~$a$. Nevertheless, in this way we can compute the Green function of the Helmholtz equation for a quite wide range of indices.\\

{\bf Acknowledgements.} We are thankful to the anonymous reviewers, whose penetrating questions and remarks have helped us to improve the manuscript.

\section{Some helpful integral formulae}\label{sec:appendix}

In this section, we collect several hints for computing integrals of irrational functions that arise when one wants to derive a formula for the kernel~$G$ in an even-dimensional space.

\begin{lemma} If $\lambda$ is a half-integer, $\lambda\in\mathbb N/2\setminus\mathbb N$, then the following holds:
\begin{align*}
\int&\left[\frac{1-r^2}{r(1-2tr+r^2)^{\lambda+1}}-\frac{1}{r}\right]dr
   =\frac{1}{\lambda(\mathbf t+\mathbf r^2)^\lambda}+\frac{1}{2}\sum_{j=1}^{\lambda-\frac{1}{2}}\frac{1}{(\lambda-j)(\mathbf t+\mathbf r^2)^{\lambda-j}}\\
&+\sum_{j=1}^{\lambda-\frac{1}{2}}\frac{t\,Q^\lambda_{j-1}(\mathbf t)\,\mathbf r}{\mathbf t^j(\mathbf t+\mathbf r^2)^{\lambda-j}}
   -\ln\left(\mathbf t-t\mathbf r+\sqrt{\mathbf t+\mathbf r^2}\right)+C,
\end{align*}
where $\mathbf r=r-t$ and $\mathbf t=1-t^2$, and polynomials $Q^\lambda_j$, $j=0,1,\dots,\lambda-\frac32$, are given recursively by
\begin{align*}
1=&\,2(\lambda-1)Q^\lambda_0(\mathbf t),\\
0=&\,2(\lambda-2)Q^\lambda_1(\mathbf t)-\left[(2\lambda-3)+2(\lambda-1)\mathbf t\right]Q^\lambda_0(\mathbf t),\\
0=&\,2(\lambda-j-1)Q^\lambda_j(\mathbf t)-\left[(2\lambda-2j-1)+2(\lambda-j)\mathbf t\right]Q^\lambda_{j-1}(\mathbf t)\\
   &+(2\lambda-2j+1)\mathbf tQ^\lambda_{j-2}(\mathbf t),\hspace{10em}j=2,\dots,\lambda-\tfrac{3}{2}.
\end{align*}\\
\end{lemma}

\begin{lemma}Suppose $\mathbf t\ne0$ and consider the integral
$$
I_{k,J+\frac12}:=\int\frac{\mathbf R^k}{(\mathbf t+\mathbf R^2)^{J+\frac12}}\,d\mathbf R.
$$
If $k=2\kappa+1$, then
$$
I_{k,J+\frac12}=\sum_{\iota=0}^\kappa\binom{\kappa}{\iota}
   \frac{(-1)^{\kappa-\iota+1}}{2(J-\iota)-1}\cdot\frac{\mathbf t^{\kappa-\iota}}{(\mathbf t+\mathbf R^2)^{J-\iota-\frac12}}+C.
$$
If $k=2\kappa$ and $\kappa<J$, then
$$
I_{k,J+\frac12}=\mathbf t^{\kappa-J}\sum_{\iota=0}^{J-\kappa-1}\binom{J-\kappa-1}{\iota}\frac{(-1)^{J-\kappa-\iota-1}}{2(J-\iota)-1}
   \cdot\left(\frac{\mathbf R^2}{\mathbf t+\mathbf R^2}\right)^{J-\iota-\frac12}+C.
$$
If $k=2\kappa$ and $\kappa\geq J$, then
$$
I_{k,J+\frac12}=\frac{\mathbf R}{(\mathbf t+\mathbf R^2)^{J-\frac12}}\cdot P_{k,J+\frac12}(\mathbf t,\mathbf R^2)
   +a\,\mathbf t^{\kappa-J}\,\ln\left(\mathbf R+\sqrt{\mathbf t+\mathbf R^2}\right)+C,
$$
where
$$
a^{\kappa,J+\frac12}=\frac{(-1)^{\kappa-J}(2\kappa-1)!!}{2^{\kappa-J}\cdot(\kappa-J)!\cdot(2J-1)!!},
$$
$P_{k,J+\frac12}$ is a homogeneous polynomial of degree~$\kappa-1$ given by
$$
P_{k,J+\frac12}(\mathbf t,\mathbf R^2)=\sum_{\iota=0}^{\kappa-1} a_\iota^{\kappa,J+\frac12}\mathbf t^{\kappa-\iota-1}\,\mathbf R^{2\iota}
$$
with coefficients obtained recurrently
\begin{align*}
a_0^{\kappa,J+\frac12}&=-a^{\kappa,J+\frac12},\\
a_1^{\kappa,J+\frac12}&=-\frac{(3J-2)\,a^{\kappa,J+\frac12}}{3},\\
a_\iota^{\kappa,J+\frac12}&=\frac{2(J-\iota)\,a_{\iota-1}^{\kappa,J+\frac12}-a^{\kappa,J+\frac12}\binom{J}{\iota}}{2\iota+1},\qquad\iota=2,3,\dots,\kappa-1,
\end{align*}
and $P_{0,J+\frac12}\equiv0$ (with the conventions $(-1)!!=1$ and $\binom{J}{\iota}=0$ for $\iota>J$).\\
\end{lemma}

\begin{lemma}\label{lem:integral_LJ12}The integral
$$
\mathcal I_{L,J+\frac12}(t,R):=\int\frac{R^L}{(1-2tR+R^2)^{J+\frac12}}\,dR,
$$
$L,J\in\mathbb N_0$, is equal to
\begin{align*}
\mathcal I_{L,J+\frac12}&=\frac{A_{L,J+\frac12}(t,R)}{(1-2tR+R^2)^{J-\frac{1}{2}}\cdot(1-t^2)^J}\\
&+B_{L,J+\frac12}(t)\cdot\ln\left(R-t+\sqrt{1-2tR+R^2}\right)+C,
\end{align*}
where the polynomials $A$ and $B$ are given by
\begin{align*}
A_{L,J+\frac12}(t,R)&=\sum_{\kappa=0}^{J-1}\sum_{\iota=0}^{J-\kappa-1}\alpha_{\kappa,\iota}^{L,J+\frac12}\,
   t^{L-2\kappa}\,\mathbf t^\kappa\,\mathbf R^{2(J-\iota-\frac12)}(\mathbf t+\mathbf R^2)^\iota\\
&+\sum_{\kappa=J}^{\left[\frac{L}{2}\right]}\sum_{\iota=0}^{\kappa-1}\beta_{\kappa,\iota}^{L,J+\frac12}\,t^{L-2\kappa}\,\mathbf t^{J+\kappa-\iota-1}\,\mathbf R^{2\iota+1}\\
&+\sum_{\kappa=0}^{\left[\frac{L-1}{2}\right]}\sum_{\iota=0}^\kappa\gamma_{\kappa,\iota}^{L,J+\frac12}\,t^{L-2\kappa-1}\,\mathbf t^{J+\kappa-\iota}\,(\mathbf t+\mathbf R^2)^\iota\\
B_{L,J+\frac12}(t)&=\sum_{\kappa=J}^{\left[\frac{L}{2}\right]}\mu_\kappa^{L,J+\frac12}\cdot t^{L-2\kappa}\cdot\mathbf t^{\kappa-J},
\end{align*}
with $\mathbf R=R-t$, $\mathbf t=1-t^2$, and
\begin{align*}
\alpha_{\kappa,\iota}^{L,J+\frac12}&=\binom{L}{2\kappa}\binom{J-\kappa-1}{\iota}\frac{(-1)^{J-\kappa-\iota-1}}{2(J-\iota)-1},\\
\beta_{\kappa,\iota}^{L,J+\frac12}&=\binom{L}{2\kappa}a_\iota^{\kappa,J+\frac12},\\
\gamma_{\kappa,\iota}^{L,J+\frac12}&=\binom{L}{2\kappa+1}\binom{\kappa}{\iota}\frac{(-1)^{\kappa-\iota+1}}{2(J-\iota)-1},\\
\mu_\kappa^{L,J+\frac12}&=\binom{L}{2\kappa}a^{\kappa,J+\frac12}.
\end{align*}
\end{lemma}

{\bf Proof.} Follows from
\begin{align*}
\mathcal I_{L,J+\frac12}&=\int\frac{(\mathbf R+t)^L}{(\mathbf t+\mathbf R^2)^{J+\frac12}}\,d\mathbf R=\sum_{k=0}^L\binom{L}{k}t^{L-k}I_{k,J+\frac12}\\
&=\sum_{\kappa=0}^{\min\left\{J,\left[\tfrac{L}{2}\right]\right\}-1}\binom{L}{2\kappa}t^{L-2\kappa}I_{2\kappa,J+\frac12}
   +\sum_{\kappa=J}^{\left[\frac{L}{2}\right]}\binom{L}{2\kappa}t^{L-2\kappa}I_{2\kappa,J+\frac12}\\
&+\sum_{\kappa=0}^{\left[\frac{L-1}{2}\right]}\binom{L}{2\kappa+1}t^{L-2\kappa-1}I_{2\kappa+1,J+\frac12}.
\end{align*}
\hfill$\Box$

\begin{lemma}
Let $k$ be a positive integer. The integral
$$
\int R^k\ln\left(1-tR+\sqrt{1-2tR+R^2}\right)\,dR
$$
equals
\begin{align*}
&p_k(t,R)\cdot\sqrt{1-2tR+R^2}+q_k(t)\cdot\ln\left(R-t+\sqrt{1-2tR+R^2}\right)\\
&+\frac{R^{k+1}}{k+1}\left[\ln\left(1-tR+\sqrt{1-2tR+R^2}\right)-\frac{1}{k+1}\right]+C,
\end{align*}
where
$$
p_k(t,R)=\sum_{j=0}^{k-1}R^j\cdot\pi_j^k(t)
$$
with $\pi_j^k$ defined recursively by
\begin{align*}
\pi_{k-1}^k(t)&=\frac{1}{k(k+1)},\\
\pi_{k-2}^k(t)&=\frac{2k-1}{k-2}\,t\pi_{k-1}^k(t),\\
\pi_j^k(t)&=\frac{1}{j+1}\left[(2j+3)t\pi_{j+1}^k(t)-(j+2)\pi_{j+2}^k(t)\right],&j=k-3,k-4,\dots1,\\
\pi_0^k(t)&=3t\pi_1^k(t)-2\pi_2^k(t),
\end{align*}
and
$$
q_k(t)=t\pi_0^k(t)-\pi_1^k(t).
$$
\end{lemma}

\end{document}